\documentclass{amsart}

\usepackage{amsmath,amscd,amssymb,amsthm}
\usepackage{latexsym,enumerate,graphicx,rotating,braket,mathrsfs,url}
\usepackage{hyperref}
\usepackage{geometry}
\usepackage[all]{xy}

\newcommand{\bibfilename}{unibib}

\numberwithin{equation}{section}

\theoremstyle{plain}
\newtheorem{thm_}[equation]{Theorem}
\newtheorem{lemma_}[equation]{Lemma}
\newtheorem{prop_}[equation]{Proposition}
\newtheorem{cor_}[equation]{Corollary}
\newtheorem{eg_}[equation]{Example}
\newtheorem{con_}[equation]{Conjecture}
\newtheorem*{cons_}{Conjecture}

\theoremstyle{definition}
\newtheorem{thmu_}[equation]{Theorem}
\newtheorem*{thmus_}{Theorem}
\newtheorem{propu_}[equation]{Proposition}
\newtheorem*{propus_}{Proposition}
\newtheorem{coru_}[equation]{Corollary}
\newtheorem{lemu_}[equation]{Lemma}
\newtheorem{egu_}[equation]{Example}
\newtheorem*{egus_}{Example}
\newtheorem{def_}[equation]{Definition}
\newtheorem*{defs_}{Definition}

\theoremstyle{remark}
\newtheorem{rk_}[equation]{Remark}

\newcommand{\thm}[1]{\begin{thm_}#1\end{thm_}}
\newcommand{\thmu}[1]{\begin{thmu_}#1\end{thmu_}}
\newcommand{\thmus}[1]{\begin{thmus_}#1\end{thmus_}}
\newcommand{\lemm}[1]{\begin{lemma_}#1\end{lemma_}}

\newcommand{\eg}[1]{\begin{eg_}#1\end{eg_}}

\newcommand{\prop}[1]{\begin{prop_}#1\end{prop_}}

\newcommand{\rk}[1]{\begin{rk_}#1\end{rk_}}

\newcommand{\pf}[1]{\begin{proof}#1\end{proof}}

\DeclareMathOperator{\Gal}{Gal}

\DeclareMathOperator{\Spec}{Spec}

\newcommand{\fracl}[3]{\genfrac{(}{)}{}{}{#1}{#2}_{#3}}
\newcommand{\fracn}[2]{\genfrac{(}{)}{}{}{#1}{#2}}

\newcommand{\ZZ}{\mathbb Z}
\newcommand{\QQ}{\mathbb Q}
\newcommand{\RR}{\mathbb R}

\renewcommand{\AA}{\mathbb A}
\newcommand{\II}{\mathbb I}
\newcommand{\GG}{\mathbb G}
\newcommand{\XX}{\mathbf X}
\newcommand{\TT}{\mathbf T}

\renewcommand{\o}{\mathfrak o}
\newcommand{\p}{\mathfrak p}
\renewcommand{\P}{\mathfrak P}

\renewcommand{\t}{\times}

\newcommand{\lr}{\longrightarrow}

\newcommand{\eq}[1]{\begin{equation}#1\end{equation}}
\newcommand{\eqn}[1]{\begin{equation*}#1\end{equation*}}

\newcommand{\aln}[1]{\begin{align*}#1\end{align*}}

\newcommand{\cs}[1]{\begin{cases}#1\end{cases}}

\newcommand{\enmt}[1]{\begin{enumerate}#1\end{enumerate}}
\renewcommand{\it}{\item}

\newcommand{\itm}[1]{\it[\upshape{(#1)}]}

\begin{document}
\title[Integral Representation of Binary Quadratic Forms]{On the Integral Representation of Binary Quadratic Forms and the Artin Condition}
\author[C. Lv]{Chang Lv}
\address{State Key Laboratory of Information Security\\
Institute of Information Engineering\\
Chinese Academy of Sciences\\
Beijing 100093, P.R. China}
\email{lvchang@amss.ac.cn}
\author[J. Shentu]{Junchao Shentu}
\address{Shanghai Center for Mathematical Science\\
Fudan University\\
Shanghai, P. R. China}
\email{stjc@amss.ac.cn}
\author[Y. Deng]{Yingpu Deng}
\address{Key Laboratory of Mathematics Mechanization\\
NCMIS, Academy of Mathematics and Systems Science\\
Chinese Academy of Sciences\\
Beijing 100190, P.R. China}
\email{dengyp@amss.ac.cn}
\date{\today}
\begin{abstract}
For diophantine equations of the form   $ax^2+bxy+cy^2+g=0$ over $\ZZ$ whose coefficients satisfy some assumptions, we show that 
a condition with respect to Artin reciprocity map, which we call the Artin condition,
is the only obstruction to the local-global principle for integral solutions of the equation. 
Some concrete examples are presented.
\end{abstract}
\maketitle

\section{Introduction}\label{sec_intro}
The main theorem of a book by  Cox~\cite{cox} is a beautiful criterion of the solvability of the diophantine equation $p=x^2+ny^2$. The specific statement is
\thmus{
Let $n$ be a positive integer. Then there is a monic irreducible polynomial $f_n(x)\in\ZZ[x]$ of degree $h(-4n)$ such that if an odd prime $p$ divides neither $n$ nor the discriminant of $f_n(x)$, then $p=x^2+ny^2$ is solvable over $\ZZ$ if and only if $\fracl{-n}{p}{}=1$ and $f_n(x)=0$ is solvable over $\ZZ/p\ZZ$. Here $h(-4n)$ is the class number of primitive positive definite binary forms of discriminant $-4n$. Furthermore, $f_n(x)$ may be taken to be the minimal polynomial of a real algebraic integer $\alpha$ for which $L=K(\alpha)$ is the ring class field of the order $\ZZ[\sqrt{-n}]$ in the imaginary quadratic field $K=\QQ(\sqrt{-n})$.
}
There are some generalizations considering the problem over quadratic fields.
By using classical results in the class field theory, the  first and third author~\cite{rcf} gave the criterion of 
the integral solvability of the  equation $p = x^2 + ny^2$  for some $n$ over a class of
imaginary quadratic fields, where $p$ is a prime element.

Recently, Harari~\cite{bmob} showed that the Brauer-Manin obstruction is the only obstruction for the existence of integral points of a scheme over the ring of integers of a number field, whose generic fiber is a principal homogeneous space (torsor) of a torus.
After then  Wei and  Xu  \cite{multi-norm-tori,multip-type}  construct the idele groups which
are the so-called $\XX$-\emph{admissible subgroups} for determining the integral points for multi-norm  tori,
and  interpret the $\XX$-admissible subgroup in terms of
finite Brauer-Manin obstruction.
In \cite[Section 3]{multi-norm-tori} Wei and Xu also showed how to apply this method to binary quadratic diophantine equations. 
As applications,  they gave  some explicit criteria of the solvability of equations of the form $x^2\pm dy^2=a$ over
$\ZZ$ in \cite[Sections 4 and 5]{multi-norm-tori}.

Later  Wei \cite{wei_diophantine} applied the method in \cite{multi-norm-tori} 
to  give some additional criteria of the solvability of the diophantine equation 
$x^2-dy^2=a$ over $\ZZ$ for some $d$.
He also determined which integers can be written as a sum of two integral squares for some of the quadratic fields $\QQ(\sqrt{\pm p})$ (in \cite{wei1}), $\QQ(\sqrt{-2p})$ (in \cite{wei2}) and so on. 

In this article, we apply the method in \cite{multi-norm-tori} to diophantine equations of the form
\eq{\label{eq_qr}
ax^2+bxy+cy^2+g=0
}
 over $\ZZ$, a binary quadratic form representing an integer.
With some additional assumptions, by choosing   $\XX$-admissible  subgroups for \eqref{eq_qr}
the same as in \cite[Sections 4, 5]{multi-norm-tori} and \cite{wei1},
 we obtain  criteria of the solvability of \eqref{eq_qr}, 
as a variant of \cite[Proposition 4.1]{multi-norm-tori} and 
\cite[Proposition 5.1]{multi-norm-tori}.
In the case $b=0$, the first and second author \cite{lv2015integral}
also gave some corresponding results.

Specifically, the main results of this article are 
(for notation one can see Section \ref{sec_nota}): 
\thmus{
Let $a,b,c$ and $g$ be integers such that $d=4ac-b^2<0$. 
Suppose $-d=p_1^{m_1}p_2^{m_2}\dots p_r^{m_r}$ where
$r>0,m_k\ge1$ are not all even and $p_k$ are distinct odd primes such that one of the 
following assumptions holds:
\enmt{
\itm{1} $p_i\equiv3\pmod4$ for some $i$.
\itm{2} $r=2$ or $r>3$ is odd, $p_i\equiv1\pmod4, m_i=1$ for all $i$ and $(p_i/p_j)=-1$ for all  $i\neq j$.
}
Set $E=\QQ(\sqrt{-d})$,
$L=\ZZ+ \ZZ\sqrt{-d} $ and $H_L$ the ring class field corresponding to $L$. 
Let $\XX=\Spec(\ZZ[x,y]/(ax^2+bxy+cy^2+g))$. 
Then $\XX(\ZZ)\neq\emptyset$ if and only if there exists
\eqn{
\prod_{p\le \infty}(x_p, y_p)\in\prod_{p\le \infty}\XX(\ZZ_p)
}
such that
\eqn{
\psi_{H_L/E}(\tilde f_E(\prod_p(x_p,y_p)))=1.
}
}
\thmus{
Let $a,b,c$ and $g$ be integers such that $d=4ac-b^2<0$. 
Suppose $-d=p_1^{m_1}p_2^{m_2}\dots p_r^{m_r}$ where
$r>0,m_k\ge1$ are not all even and $p_k$ are distinct odd primes such that one of the 
following assumptions holds:
\enmt{
\itm{1} $p_i\equiv3\pmod4$ for some $i$.
\itm{2} $r=2$ or $r>3$ is odd, $p_i\equiv1\pmod4, m_i=1$ for all $i$ and $(p_i/p_j)=-1$ for all  $i\neq j$.
}
Set $E=\QQ(\sqrt{-d})$,
$L=\ZZ+ \ZZ\sqrt{-d} $ and $H_L$ the ring class field corresponding to $L$. 
Let $\XX=\Spec(\ZZ[x,y]/(ax^2+bxy+cy^2+g))$. 
Then $\XX(\ZZ)\neq\emptyset$ if and only if there exists
\eqn{
\prod_{p\le \infty}(x_p, y_p)\in\prod_{p\le \infty}\XX(\ZZ_p)
}
such that
\eqn{
\psi_{H_L/E}(\tilde f_E(\prod_p(x_p,y_p)))=1.
}
}

In Section \ref{sec_artin_cond}, we introduce  from  \cite{multi-norm-tori}  notation
 and the  general result we mainly use in this paper, but in a modified way which focus on our goal.
Then we give our results on the equation  \eqref{eq_qr} in Section \ref{sec_z}.
If the discriminant $d$ is positive we need  no additional assumption.
But  if $d$ is negative, we add some assumptions on it
(as it is done in \cite[Section 5]{multi-norm-tori}).
The results state that
the integral local solvability and the Artin condition
(see Remark \ref{rk_artin_cond})
 completely describe the global integral solvability.
We also give some examples  showing  the explicit criteria of the solvability.

\section{Solvability by the Artin Condition}\label{sec_artin_cond}

\subsection{Notation}\label{sec_nota}
Let $F$ be a number field, $\o_F$ the ring of integers of $F$,
$\Omega_F$ the set of all places in $F$ and $\infty_F$ the set of all infinite places in $F$. 
Let $F_\p$ be the completion of $F$ at $\p$ and $\o_{F_\p}$ be the valuation ring of $F_\p$ for each $\p\in\Omega_F\setminus\infty_F$. We also write $\o_{F_\p}=F_\p$ for $\p\in\infty_F$.
The adele ring (resp.  idele group) of $F$ is denoted by $\AA_F$ (resp. $\II_F$).

Let $a,b,c$ and $g$ be  elements in $\o_F$ and suppose that $-d=b^2-4ac$ is not a square in $F$. Let $E=F(\sqrt{-d})$ and
 $\XX=\Spec(\o_F[x,y]/(ax^2+bxy+cy^2+g))$ be the affine scheme defined by the equation $ax^2+bxy+cy^2+g=0$ over $\o_F$.
The equation 
\eq{\label{eq_bqf}
ax^2+bxy+cy^2+g=0
}
is solvable over $\o_F$ if and only if $\XX(\o_F)\neq\emptyset$.

Now we denote 
\aln{
\tilde x&:= 2ax+by,	\\
\tilde y&:= y,	\\
n&:=-4ag.
}
Then we can write  \eqref{eq_bqf} as
\eq{\label{eq_bqf_n}
\tilde x^2+d\tilde y^2=n.
}

Denote $R_{E/F}(\GG_m)$ the Weil restriction of $\GG_{m,E}$ to $F$. Let
\eqn{ \varphi: R_{E/F}(\GG_m)\lr\GG_m }
be the homomorphism of algebraic groups which represents
\eqn{ x\mapsto N_{E/F}(x): (E\otimes_FA)^\t\lr A^\t}
for any $F$-algebra $A$. Define the torus $T:=\ker\varphi$. 
Let $X_F$ be the generic fiber of $\XX$.
We can identify elements in $T(A)$ (resp. $(X_F(A)$) as
$u+\sqrt{-d}v$ (resp. $\tilde x+\sqrt{-d}\tilde y$).
Then $X_F$ is naturally a $T$-torsor by the action: 
\aln{
T(A)\t X_F(A) &\lr X_F(A)\\
(u+\sqrt{-d}v, \tilde x+\sqrt{-d}\tilde y) &\mapsto (u+\sqrt{-d}v) (\tilde x+\sqrt{-d}\tilde y).
}
Note that $T$ has an integral model $\TT=\Spec(\o_F[x,y]/(x^2+dy^2-1))$
and we can view $\TT(o_{F_\p})$ as a subgroup of $T(F_\p)$.

Denote by $\lambda$ the embedding of $T$ into $R_{E/F}(\GG_m)$.
Clearly $\lambda$ induces a natural injective group homomorphism
\eqn{ \lambda_E: T(\AA_F)\lr\II_E. }
Let $L=\o_F+\o_F\sqrt{-d}$ in $E$ and $L_\p=L\otimes_{\o_F}\o_{F_\p}$
in  $E_\p=E\otimes_F F_\p$.
Then
\eqn{ \TT(\o_{F_\p})=\set{\beta\in L_\p^\t|N_{E_\p/F_\p}(\beta)=1}.}
It follows that  $\lambda_E(\TT(\o_{F_\p})) \subseteq  L_\p^\t.$
Note that $\lambda_E( T(F)) \subseteq E^\t$ in $\II_E$. 
Let $\Xi_L:=\prod_{\p\in\Omega_F} L_\p^\t$ which is an open  subgroup of $\II_E$.
Then  the following map  induced by $\lambda_E$ is well-defined:
\eqn{
\tilde\lambda_E: T(\AA_F)/T(F)\prod_{\p\in \Omega_F}\TT( \o_{F_\p})\lr \II_E/E^\t \prod_{\p\in\Omega_F} L_\p^\t.
}
Now we assume that 
\eq{\label{eq_nonempty}
\XX(F)\neq\emptyset,
}
i.e. $X_F$ is a trivial $T$-torsor.
Fixing a rational point $P\in X_F(F)$,  
for any $F$-algebra $A$, we have an isomorphism
\aln{
\phi_P: X_F(A) &\cong T(A)\\
x &\mapsto P^{-1}x
}
induced by $P$. 
Since we can view $\prod_{\p\in\Omega_F}\XX(\o_{F_\p})$ as a subset of $X_F(\AA_F)$, the composition $f_E:=\lambda_E\phi_P: \prod_\p\XX(\o_{F_\p})\lr  \II_E$ makes sense, mapping 
$x$ to $P^{-1}x$ in $\II_E$. Note that $P$ is in $E^\t\subset \II_E$ since it is a rational point over $F$. It follows that we can 
define the map $\tilde f_E$ to be the composition
\eqn{\xymatrix{
\prod_\p\XX(\o_{F_\p}) \ar@{->}[r]^-{ f_E}  &\II_E       \ar@{->}[r]^-{\t P} &\II_E\\
x                                 \ar@{|->}[r]              &P^{-1}x     \ar@{|->}[r]       &x.
}}
It can be seen that the restriction to $\XX(\o_{F_\p})$ of   $\tilde f_E$ is defined by
\eq{\label{eq_tilde_f_E}
\tilde f_E[(x_\p,y_\p)]=
\cs{
(\tilde x_\p+\sqrt{-d}\tilde y_\p,   \tilde x_\p-\sqrt{-d}\tilde y_\p) \in E_{\P_1}\t E_{\P_2}	&\text{if }\p\text{ splits in }E/F,\\
\tilde x_\p+\sqrt{-d}\tilde y_\p \in E_\P										&\text{otherwise},
}}
where $\P_1$ and $\P_2$ (resp. $\P$) are places of $E$ above $\p$.

Recall that $L=\o_F+\o_F\sqrt{-d}$, $L_\p=L\otimes_{\o_F}\o_{F_\p}$ and  $\Xi_L=\prod_\p L_\p^\t$ is an open subgroup of $\II_E$.  By the \emph{ring class field corresponding to $L$} we mean the class field $H_L$ corresponding to $\Xi_L$ under the class field theory,
such that  the Artin map  gives the isomorphism $\psi_{H_L/E}: \II_E/E^\t\Xi_L\cong\Gal(H_L/E)$.
For any $\prod_{\p\in\Omega_F}(x_\p, y_\p)\in\prod_{\p\in\Omega_F}\XX(\o_{F_\p})$, noting that $P$ is in $E$, we have 
\eq{\label{eq_tilde_f}
\psi_{H_L/E}(f_E(\prod_\p(x_\p,y_\p)))=1\text{ if and only if }
\psi_{H_L/E}(\tilde f_E(\prod_\p(x_\p,y_\p)))=1.
}

\rk{\label{rk_hasse_min}
If  $\prod_{\p\in\Omega_F}\XX(\o_{F_\p})\neq\emptyset$,  then the assumption \eqref{eq_nonempty} 
holds automatically  by the Hasse-Minkowski theorem on quadratic equations. 
Hence We can pick an $F$-point $P$ of $X_F$ and obtain $\phi_P$.
}

\subsection{A general result}\label{sec_main}
In the previous section, we choose the subgroup to be  $\Xi_L=\prod_\p L_\p^\t$ 
where $L=\o_F+\o_F\sqrt{-d}$  and $L_\p=L\otimes_{\o_F}\o_{F_\p}$. 
By some additional assumptions, we  prove that $\Xi_L$ can be viewed as an
admissible subgroup for    $\XX=\Spec(\o_F[x,y]/(ax^2+bxy+cy^2+g))$. 
\newcommand{\U}{\mathcal U}
\lemm{\label{lambda_inj}
Let $\U$ be a complete system of representatives of $\o_F^\t/(\o_F^\t)^2$.
Suppose  for every $u\in\U$, 
the equation $x^2+dy^2=u$ is solvable over $\o_F$
or is not solvable over $\o_{F_\p}$ for some place $\p$.
Then the map $\tilde\lambda_E$ is injective.
}
\pf{
Recall that $T=\ker(R_{E/F}(\GG_m)\lr\GG_m)$ and $\TT$ is the group scheme 
defined by the equation $x^2+dy^2=1$ over $\o_F$.
  Therefore we have
\eqn{ T(F)=\set{\beta\in E^\t|N_{E/F}(\beta)=1} }
and 
\eqn{ \TT(\o_{F_\p})=\set{\beta\in L_\p^\t|N_{E_\p/F_\p}(\beta)=1}. }
Suppose   $t\in T(\AA_F)$ such that $\tilde\lambda_E(t)=1$. Write $t=\beta i$ with $\beta\in E^\t$ and $i\in \prod_\p L_\p^\t$.
Since $t\in T(\AA_F)$ we have
\eqn{ N_{E/F}(\beta)N_{E/F}(i)=N_{E/F}(\beta i)=1. }
It follows that
\eqn{ N_{E/F}(i)=N_{E/F}(\beta^{-1})\in F^\t\cap \prod_\p\o_{F_\p}^\t=\o_F^\t. }
So by the definition of $\U$, we have 
  $N_{E/F}(i)=uv^2$ for some $u\in\U$ and $v\in\o_F^\t$.
Then 
\eqn{
N_{E/F}(iv^{-1})=u,
}
from which we know 
that  the equation $x^2+dy^2=u$ is  solvable over $\o_{F_\p}$ for every place $\p$ of $F$,
since  $v^{-1}\in\o_F$.
Thus the assumption tells us  that  $x^2+dy^2=u$ is solvable over $\o_F$.
Let  $(x_0,y_0)\in \o_F^2$ be such a solution and let
\aln{ 
\zeta&=x_0+y_0\sqrt{-d}, \\
\gamma&=\beta v\zeta,\\
\text{and }j&=iv^{-1}\zeta^{-1}. 
}
Then $N_{E/F}(\gamma)=N_{E/F}(j)=1$. Note that $\zeta\in L_\o^\t$,
and we have $\gamma\in T(F)$ and $j\in\prod_\p\TT(\o_{F_\p})$. 
It follows that $t=\beta i=\gamma j\in T(F)\prod_\p\TT(\o_{F_\p})$. 
This finishes the proof.
}

As a result, we can obtain 
 criteria of the solvability in a more explicit way. 
We state it in the following proposition, which is a Corollary to
\cite[Corollary 1.6]{multi-norm-tori}.
\prop{\label{prop_artin_cond}
Let symbols be as before and $\U$ satisfy the assumption in 
Lemma \ref{lambda_inj}.
Then $\XX(\o_F)\neq\emptyset$ if and only if there exists
\eqn{
\prod_{\p\in\Omega_F}(x_\p, y_\p)\in\prod_{\p\in\Omega_F}\XX(\o_{F_\p})
}
such that
\eq{\label{eq_artin_cond}
\psi_{H_L/E}(\tilde f_E(\prod_\p(x_\p,y_\p)))=1.
}}
\pf{
By the assumption 
we know from Lemma \ref{lambda_inj}, that
\eq{ \label{eq_lambda_inj}
\tilde \lambda_E: T(\AA_F)/T(F)\prod_\p \TT(\o_{F_\p})\lr \II_E/E^{\t}\prod_\p L_\p^\t
}
is injective.

If $\XX(\o_F)\neq\emptyset$, then 
\eqn{
\tilde f_E\left(\prod_\p\XX(\o_{F_\p})\right)\cap E^\t\prod_\p L_\p^\t\supseteq \tilde f_E(\XX(\o_F))\cap E^\t\neq\emptyset
}
Hence there exists $x\in \prod_{\p\in\Omega_F}\XX(\o_{F_\p})$ such that $\psi_{H_L/E}\tilde f_E(x)=1$.

Conversely, 
suppose there exists $x\in \prod_\p\XX(\o_{F_\p})$ such that $\psi_{H_L/E}\tilde f_E(x)=1$ 
(here $\tilde f_E$ makes sense by Remark \ref{rk_hasse_min}), i.e. 
$\lambda_E\phi_P(x)=f_E(x)\in \Xi_L=E^{\t}\prod_\p L_\p^\t$.
Since $\tilde \lambda_E$, i.e. \eqref{eq_lambda_inj}, is injective, there are $\tau\in T(F)$ and $\sigma\in\prod_\p\TT( \o_{F_\p})$ such that $\tau\sigma=\phi_P(x)=P^{-1}x$,
i.e. $\tau\sigma(P)=x$. Since $P\in X_F(F)$ and 
\eq{\label{eq_in_stab}
g\XX(\o_{F_\p})=\XX(\o_{F_\p})\text{ for all }g\in\TT(\o_{F_\p}),
}
 it follows that 
\eqn{
\tau(P)=\sigma^{-1}(x)\in \XX(F)\cap \prod_\p\XX(\o_{F_\p})=\XX(\o_F).
}
Then the proof is done.
}
\rk{\label{rk_artin_cond}
The condition \eqref{eq_artin_cond} is called the \emph{Artin condition} in, for example, Wei's~\cite{wei_diophantine,wei1,wei2}. It interprets the fact that the Brauer-Manin obstruction is the only obstruction for existence of the integral points by conditions in terms of the class field theory.
Consequently, if the assumption  in the proposition holds, 
 the integral local solvability  and the Artin condition completely describe the global integral solvability.
As a result, in cases where the ring class fields are known it is possible to calculate the Artin condition, giving  explicit criteria for the solvability.
}

\section{The Integral Representation of Binary Quadratic Forms over $\ZZ$} \label{sec_z}
Now we consider the case where $F=\QQ$ which is our focus. We now distinguish  the sign of the discriminant $d$.

\subsection{The case where the discriminant  $d>0$}
\thm{\label{thm_q}
Let $a,b,c$ and $g$ be integers and suppose that $d=4ac-b^2>0$. Set $E=\QQ(\sqrt{-d})$,
$L=\ZZ+ \ZZ\sqrt{-d} $ and $H_L$ the ring class field corresponding to $L$. 
Let $\XX=\Spec(\ZZ[x,y]/(ax^2+bxy+cy^2+g))$. 
Then $\XX(\ZZ)\neq\emptyset$ if and only if there exists
\eqn{
\prod_{p\le \infty}(x_p, y_p)\in\prod_{p\le \infty}\XX(\ZZ_p)
}
such that
\eqn{
\psi_{H_L/E}(\tilde f_E(\prod_p(x_p,y_p)))=1.
}
}
\pf{
Since $d>0$ it is clear that $x^2+dy^2=-1$ is not solvable over $\RR$, which is to 
say the assumption in Proposition \ref{prop_artin_cond} holds since the only units of $\ZZ$ are $\{\pm1\}$. 
Then the result follows from  Proposition \ref{prop_artin_cond}.
}
We now give an example where the explicit  criterion is obtained using this result.
\eg{
Let $g$ be a negative integer and $l(x)=x^4-x^3+x+1\in\ZZ[x]$. 
Write $g=-2^{s_1}\t7^{s_2}\t\prod_{k=1}^r p_k^{m_k}$ ,
where $s_1,s_2,k\ge0,m_k\ge1,p_1,p_2,\dots,p_r \neq 2,7$ are distinct primes.
Define $C=\set{3, p_1, p_2, \dots, p_r}$ and
\eqn{
D=\set{p\in C | \fracn{-14}{p}=1\text{ and }l(x)\mod p\text{ irreducible}}.
}
Then the diophantine equation $3x^2+2xy+5y^2+g=0$ is solvable over $\ZZ$ if and only if
\enmt{
\itm{1} $3g\t 2^{-s_1}\equiv\pm1\pmod8$,
\itm{2} $\fracn{g\t 7^{-s_2}}{7} = 1$,
\itm{3} for all $p\nmid2\t3\t7$ with odd $m_p:= v_p(g)$, $\fracn{-14}{p}=1$,
\itm{4} and $\sum_{p\in D}v_p(3g)\equiv0\pmod2$.
}}
\pf{
In this example, we have $a=3,b=2,c=5,d=4ac-b^2=4\t14$.  Let $E= \QQ(\sqrt{-d})$.
Since $b=2$, we can simplify the equation \eqref{eq_bqf_n} by canceling $4$ in both sides.
Thus we set
\aln{
n&=-4ag/4=-3g,\\
\tilde x&=(2ax+by)/2= 3x+y,\\
    \tilde y&=y.
}
In fact  we may assume $d=14$ and Theorem \ref{thm_q} still applies.
Because if $d=14$,  we still have $E = \QQ(\sqrt{-d})$, $\tilde x^2+d\tilde y^2=n$ and also \eqref{eq_in_stab} holds.
It follows that  $L=\ZZ+\ZZ\sqrt{-14}=\o_E$ and $H_L=H_E=E(\alpha)$ the Hilbert field of $E$ 
where the minimal polynomial of $\alpha$ is $l(x)$. The Galois group \eqn{ \Gal(H_L/E)=<\sqrt{-1}>\cong\ZZ/4\ZZ.}
Let $\XX=\Spec(\ZZ[x,y]/(3x^2+2xy+5y^2+g))$ and 
\eqn{
\tilde f_E[(x_p,y_p)]=
\cs{
(\tilde x_p+\sqrt{-14}\tilde y_p,   \tilde x_p-\sqrt{-14}\tilde y_p)	&\text{if }p\text{ splits in }E/\QQ,\\
 \tilde x_p+\sqrt{-14}\tilde y_p                      &\text{otherwise},
}}
Then by Theorem \ref {thm_q}, $\XX(\ZZ)\neq\emptyset$ if and only if there exists
\eqn{
\prod_{p\le \infty}(x_p, y_p)\in\prod_{p\le \infty}\XX(\ZZ_p)
}
such that
\eqn{
\psi_{H_L/E}(\tilde f_E(\prod_p(x_p,y_p)))=1.
}
Next we calculate these conditions in details.
Recall that $n=-3g$.
By  a simple  calculation we know the local condition  \eqn{
\prod_{p\le \infty}\XX(\ZZ_p)\neq\emptyset
} is equivalent to 
\eq{\label{eq_eg_q_local}
\cs{
n\t 2^{-s_1}\equiv\pm1\pmod8,\\
\fracn{n\t 7^{-s_2}}{7} = 1,\\
\text{for all }p\nmid2\t3\t7\text{ with odd $m_p=v_p(n)$},  \fracn{-14}{p}=1.
}}
For the Artin condition, let $(x_p,y_p)_p\in\prod_p\XX(\ZZ_p)$. Then 
\eq{\label{eq_local_decomp}
(\tilde x_p+\sqrt{-14}\tilde y_p)(  \tilde x_p-\sqrt{-14}\tilde y_p)=  n\text{ in }E_\P\text{ with }\P\mid p,
}
and since $H_L/E$ is unramified, for any $p\neq \infty$ we have 
\eq{\label{eq_psi_1}
1=\cs{
   \psi_{H_L/E}(p_\P)\psi_{H_L/E}(p_{\bar\P}),    &\text{ if }p=\P\bar\P\text{ splits in }E/\QQ,\\
   \psi_{H_L/E}(p_\P),                             &\text{ if }p=\P\text{ is inert in }E/\QQ,
}}
where $p_\P$ (resp. $p_{\bar\P}$)  is in $\II_E$ such that its $\P$ (resp. $\bar\P$) component is $p$ and the other
components are $1$.
We calculate $\psi_{H_L/E}(\tilde f_E[(x_p,y_p)])$ separately:
\enmt{
\it If $p=2$, $2=\P_2^2$ in $E/\QQ$. Suppose $\P_2=\pi_2\o_{E_{\P_2}}$ for $\pi_2\in\o_{E_{\P_2}}$.
Noting that $H_L/E$ is unramified,
since $\P_2^2$ is principal in $E$ but $\P_2$ is not, we have $\psi_{H_L/E}((\pi_2)_{\P_2})=-1$.
By \eqref{eq_local_decomp} we have
\eqn{ v_{\P_2}(\tilde x_2+\sqrt{-14}\tilde y_2)=v_{\P_2}(\tilde x_2-\sqrt{-14}\tilde y_2)=\frac{1}{2}v_{\P_2}(n)=v_2(n)=s_1.}
It follows that \aln{ 
\psi_{H_L/E}(\tilde f_E[(x_2,y_2)])&=\psi_{H_L/E}((\tilde x_2+\sqrt{-14}\tilde y_2)_{\P_2})\\
                                   &=(-1)^{v_{\P_2}(\tilde x_2+\sqrt{-14}\tilde y_2)}=(-1)^{s_1},
}where $\tilde f_E[(x_2,y_2)]$ is also regarded as an element in $\II_E$ such that the component above $2$ is given by the value 
of $\tilde f_E[(x_2,y_2)]$ and $1$ otherwise.
\it If $p=7$, a similar argument shows that $\psi_{H_L/E}(\tilde f_E[(x_7,y_7)])=(-1)^{s_2}$.
\it If $\fracn{-14}{p}=1$  then by \eqref{eq_psi_1} we can distinguish the following cases:
\enmt{[(i)]
    \it $l(x)\mod p$ splits into linear factors.
    Then $\psi_{H_L/E}(p_\P)=\psi_{H_L/E}(p_{\bar\P})=1$ and  $\psi_{H_L/E}(\tilde f_E[(x_p,y_p)])=1$.
    \it $l(x)\mod p$ splits into two irreducible factors. 
    Then $\psi_{H_L/E}(p_\P)=\psi_{H_L/E}(p_{\bar\P})=-1$.
    It follows that \aln{ 
    \psi_{H_L/E}(\tilde f_E[(x_p,y_p)])&=\psi_{H_L/E}((\tilde x_p+\sqrt{-14}\tilde y_p)_\P)\psi_{H_L/E}((\tilde x_p-\sqrt{-14}\tilde y_p)_{\bar\P})\\
                                       &=(-1)^{v_\P(\tilde x_p+\sqrt{-14}\tilde y_p)+v_{\bar\P}(\tilde x_p-\sqrt{-14}\tilde y_p)}=(-1)^m,
    }where $m= v_p(n)$ since \aln{
    v_\P(\tilde x_p+\sqrt{-14}\tilde y_p)&+v_{\bar\P}(\tilde x_p-\sqrt{-14}\tilde y_p)\\
                            &=v_p(\tilde x_p+\sqrt{-14}\tilde y_p)+v_p(\tilde x_p-\sqrt{-14}\tilde y_p)=v_p(n).
    }
    \it $l(x)\mod p$  is  irreducible.
    Then $\psi_{H_L/E}(p_\P)=-\psi_{H_L/E}(p_{\bar\P})=\pm\sqrt{-1}$.
    It follows that \aln{ 
    \psi_{H_L/E}&(\tilde f_E[(x_p,y_p)])=\psi_{H_L/E}((\tilde x_p+\sqrt{-14}\tilde y_p)_\P)\psi_{H_L/E}((\tilde x_p-\sqrt{-14}\tilde y_p)_{\bar\P})\\
                &=(\pm\sqrt{-1})^{v_\P(\tilde x_p+\sqrt{-14}\tilde y_p)+v_{\bar\P}(\tilde x_p-\sqrt{-14}\tilde y_p)}
                (-1)^{v_{\bar\P}(\tilde x_p-\sqrt{-14}\tilde y_p)}\\
                &=(\pm\sqrt{-1})^m(-1)^u
    }where $m =v_p(n)$  and $u=v_p(\tilde x_p-\sqrt{-14}\tilde y_p)$ (in $\QQ_p$, $0\le u\le m$). 
By Hensel's lemma, we can choose a local solution
    $(x_p,y_p)$ suitably, such that $u$  riches any value between $0$ and $m$. 
    Hence $ \psi_{H_L/E}(\tilde f_E[(x_p,y_p)]) =\pm(\sqrt{-1})^m$ with the sign chosen freely.
}
\it If $\fracn{-14}{p}=-1$ then $p$ is inert in $E/\QQ$. By \eqref{eq_psi_1} we have $\psi_{H_L/E}(\tilde f_E[(x_p,y_p)])=1$.
\it At last if $p=\infty$, since $H_L/E$ is unramified, we have $\psi_{H_L/E}(\tilde f_E[(x_\infty,y_\infty)])=1$.
}
Putting the above argument together,
and noting that $D \neq \emptyset$ since $3\in D$ and that $n=-3g$,
  we know the Artin condition is 
\eq{\label{eq_eg_q_artin}
\sum_{p\in D}v_p(3g)\equiv0\pmod2.
}
The proof is done if we put the local condition \eqref{eq_eg_q_local} and the Artin condition \eqref{eq_eg_q_artin} together.
}

\subsection{The case where the discriminant $d<0$}
In this case,  $x^2+dy^2=-1$ is  solvable over $\RR$, so we must look for other place $p$ of $\QQ$ such that
 $x^2+dy^2=-1$  is not solvable over $\ZZ_p$.
For a rational prime $p$ that divides $d$, 
we  observe that, by Hensel's Lemma, $x^2+dy^2=-1$  is solvable over $\ZZ_p$ if and only if
it is solvable over $\ZZ/p\ZZ$, i.e. $\fracn{-1}{p}=1$.
So if $d$ is divisible by some rational prime $p$ where $p\equiv3\pmod4$ then
 $x^2+dy^2=-1$  is not solvable over $\ZZ_p$.
Otherwise if none of the  prime divisors of $d$ are congruent to  $3$ modulo $4$, 
we hope that  $x^2+dy^2=-1$  is  solvable over $\ZZ$, in order to make the 
assumption true in Proposition \ref{prop_artin_cond}.
We need the following result by  Morris Newman \cite{morris-pell-equation}.
\thmu{\label{thm_morris_pell}
Let  $r>1$ be $2$ or  odd, $p_1,p_2,\dots,p_r$ be  distinct primes
such that 
\aln{
&p_i\equiv1\pmod4,\quad 1\le i\le r,\\
&\fracn{p_i}{p_j}=-1,\quad 1\le i\neq j\le r.
}
Then the diophantine equation $x^2-p_1p_2\dots p_r y^2=-1$ has a solution in $\ZZ$.
}
Now we have the criterion for certain $d<0$.
\thm{\label{thm_q_d_neg}
Let $a,b,c$ and $g$ be integers such that $d=4ac-b^2<0$. 
Suppose $-d=p_1^{m_1}p_2^{m_2}\dots p_r^{m_r}$ where
$r>0,m_k\ge1$ are not all even and $p_k$ are distinct odd primes such that one of the 
following assumptions holds:
\enmt{
\itm{1} $p_i\equiv3\pmod4$ for some $i$.
\itm{2} $r=2$ or $r>3$ is odd, $p_i\equiv1\pmod4, m_i=1$ for all $i$ and $(p_i/p_j)=-1$ for all  $i\neq j$.
}
Set $E=\QQ(\sqrt{-d})$,
$L=\ZZ+ \ZZ\sqrt{-d} $ and $H_L$ the ring class field corresponding to $L$. 
Let $\XX=\Spec(\ZZ[x,y]/(ax^2+bxy+cy^2+g))$. 
Then $\XX(\ZZ)\neq\emptyset$ if and only if there exists
\eqn{
\prod_{p\le \infty}(x_p, y_p)\in\prod_{p\le \infty}\XX(\ZZ_p)
}
such that
\eqn{
\psi_{H_L/E}(\tilde f_E(\prod_p(x_p,y_p)))=1.
}
}
\pf{
The  units of $\ZZ$ are $\{\pm1\}$ so we only need to consider the unit $-1$.
If $(1)$ holds, i.e. $p_i\equiv3\pmod4$ for some $i$, 
one can see immediately that $x^2+dy^2=-1$ is not solvable over $\ZZ_{p_i}$.
Otherwise $(2)$ holds and then $x^2+dy^2=-1$ is  solvable over $\ZZ$ by Theorem \ref{thm_morris_pell}.
Hence the assumption in Proposition \ref{prop_artin_cond} holds and 
we complete the proof by Proposition \ref{prop_artin_cond}.
}
We now give an example for  this case.
\eg{
Let $g$ be a nonzero integer and $l(x)=x^3-x^2-4x+2\in\ZZ[x]$. 
Write $g=\pm 2^{s_1}\t79^{s_2}\t\prod_{k=1}^r p_k^{m_k}$ ,
where $s_1,s_2,k\ge0,m_k\ge1,p_1,p_2,\dots,p_r \neq 2,79$ are distinct primes.
Define $C=\set{5,p_1,p_2,\dots,p_r}$ and 
\eqn{
D=\set{p\in C | \fracn{79}{p}=1\text{ and }l(x)\mod p\text{ irreducible}}.
}
Then the diophantine equation $5x^2+14xy-6y^2+g=0$ is solvable over $\ZZ$ if and only if
\enmt{
\itm{1} $\fracn{g\t (-79)^{-s_2}}{79} = -1$,
\itm{2} for all $p\nmid2\t5\t79$ with odd $m_p:= v_p(g)$, $\fracn{79}{p}=1$,
\itm{3} and if $\set{p\in D | v_p(5g)=1}\neq\emptyset$ then $r>1$.
}}
\pf{
In this example, we have $a=5,b=14,c=-6,d=4ac-b^2=-4\t79$.  Let $E= \QQ(\sqrt{-d})$.
Since $2\mid b$, we may  cancel $4$ in both sides and assume $d=-79$ as we do in the previous example.
Since $79\equiv 3\pmod4$ the assumption $(1)$ in  Theorem \ref{thm_q_d_neg} is correct.
It follows that we can apply the theorem for $d=-79$.
Thus we set
\aln{
n&=-4ag/4=-5g,\\
\tilde x&=(2ax+by)/2= 5x+7y,\\
    \tilde y&=y.
}
Now   $E = \QQ(\sqrt{79})$, $\tilde x^2-79\tilde y^2=n$ and
$L=\ZZ+\ZZ\sqrt{79}=\o_E$ and $H_L=H_E=E(\alpha)$ the Hilbert field of $E$ 
where the minimal polynomial of $\alpha$ is $l(x)$. The Galois group \eqn{ \Gal(H_L/E)=<\omega>\cong\ZZ/3\ZZ.}
Let $\XX=\Spec(\ZZ[x,y]/(5x^2+14xy-6y^2+g))$ and 
\eqn{
\tilde f_E[(x_p,y_p)]=
\cs{
(\tilde x_p+\sqrt{79}\tilde y_p,   \tilde x_p-\sqrt{79}\tilde y_p)	&\text{if }p\text{ splits in }E/\QQ,\\
 \tilde x_p+\sqrt{79}\tilde y_p                      &\text{otherwise},
}}
Then by Theorem \ref {thm_q}, $\XX(\ZZ)\neq\emptyset$ if and only if there exists
\eqn{
\prod_{p\le \infty}(x_p, y_p)\in\prod_{p\le \infty}\XX(\ZZ_p)
}
such that
\eqn{
\psi_{H_L/E}(\tilde f_E(\prod_p(x_p,y_p)))=1.
}
By  a simple  computation the local condition  \eqn{
\prod_{p\le \infty}\XX(\ZZ_p)\neq\emptyset
} is equivalent to the first two condition $(1)$ and $(2)$ above.
For the Artin condition, let $(x_p,y_p)_p\in\prod_p\XX(\ZZ_p)$. Then 
\eqn{
(\tilde x_p+\sqrt{79}\tilde y_p)(  \tilde x_p-\sqrt{79}\tilde y_p)=  n\text{ in }E_\P\text{ with }\P\mid p,
}
and since $H_L/E$ is unramified, for any $p\neq \infty$ we have 
\eq{\label{eq_psi_1_1}
1=\cs{
   \psi_{H_L/E}(p_\P)\psi_{H_L/E}(p_{\bar\P}),    &\text{ if }p=\P\bar\P\text{ splits in }E/\QQ,\\
   \psi_{H_L/E}(p_\P),                             &\text{ if }p=\P\text{ is inert in }E/\QQ.
}}
We calculate $\psi_{H_L/E}(\tilde f_E[(x_p,y_p)])$ separately:
\enmt{
\it If $p=2$, $2=\P_2^2$ in $E/\QQ$. Suppose $\P_2=\pi_2\o_{E_{\P_2}}$ for $\pi_2\in\o_{E_{\P_2}}$.
Noting that $H_L/E$ is unramified,
since $\P_2$ is principal in $E$, we have $\psi_{H_L/E}((\pi_2)_{\P_2})=1$.
Hence $\psi_{H_L/E}(\tilde f_E[(x_2,y_2)]) = 1$.
\it If $p=79$, a similar argument shows that $\psi_{H_L/E}(\tilde f_E[(x_{79},y_{79})])=1$.
\it If $\fracn{79}{p}=1$  then by \eqref{eq_psi_1_1} we can distinguish the following two cases:
\enmt{[(i)]
    \it $l(x)\mod p$ splits into linear factors.
    Then $\psi_{H_L/E}(p_\P)=\psi_{H_L/E}(p_{\bar\P})=1$ and  $\psi_{H_L/E}(\tilde f_E[(x_p,y_p)])=1$.
   \it $l(x)\mod p$ is  irreducible.
    Then $\psi_{H_L/E}(p_\P)=(\psi_{H_L/E}(p_{\bar\P}))^{-1}= \omega^{\pm1}$.
    It follows that \aln{ 
    \psi_{H_L/E}&(\tilde f_E[(x_p,y_p)])=\psi_{H_L/E}((\tilde x_p+\sqrt{79}\tilde y_p)_\P)\psi_{H_L/E}((\tilde x_p-\sqrt{79}\tilde y_p)_{\bar\P})\\
                &=\omega^{\pm (v_\P(\tilde x_p+\sqrt{79}\tilde y_p)+v_{\bar\P}(\tilde x_p-\sqrt{79}\tilde y_p))}
                /\omega^{\pm2 v_{\bar\P}(\tilde x_p-\sqrt{79}\tilde y_p)}\\
                &=\omega^{\pm(m-2u)}
    }where $m =v_p(n)$  and $u=v_p(\tilde x_p-\sqrt{79}\tilde y_p)$ (in $\QQ_p$, $0\le u\le m$). 
By Hensel's lemma, we can choose a local solution
    $(x_p,y_p)$ suitably, such that $u$  riches any value between $0$ and $m$. 
    Hence \eqn{
 \psi_{H_L/E}(\tilde f_E[(x_p,y_p)]) =
\cs{
1 &\text{ if } m=0,\\
\omega^{\pm1} &\text{ if } m=1,\\
1 \text{ or } \omega^{\pm1} &\text{ if } m\ge2,
}} where the values are chosen freely in each case.
}
\it If $\fracn{79}{p}=-1$ then $p$ is inert in $E/\QQ$. By \eqref{eq_psi_1_1} we have $\psi_{H_L/E}(\tilde f_E[(x_p,y_p)])=1$.
\it At last if $p=\infty$, since $H_L/E$ is unramified, we have $\psi_{H_L/E}(\tilde f_E[(x_\infty,y_\infty)])=1$.
}
Putting the above argument together,
and noting that $5\in D$ and $n=-5g$,
  we know the Artin condition is  exactly the last condition $(3)$ in the example.
This completes the proof.
}

\section*{Acknowledgment}
The author would like to thank Yupeng Jiang and Jianing Li for  helpful discussions
and the referees for valuable suggestions.

\bibliography{\bibfilename}
\bibliographystyle{amsplain} 
\end{document}